\documentclass{amsart} 

\usepackage[latin1]{inputenc}
\usepackage{amsmath}

\usepackage[round]{natbib}

\usepackage[all]{xy}

\newcommand{\cut}[1]{\overline{#1}} 
\newcommand{\fan}{\mathfrak{F}}
\newcommand{\FF}{F}
\newcommand{\NN}{\mathbb{N}}

\newcommand{\RK}{\Delta'}
\newcommand{\SK}{\Delta}

\newcommand{\RR}{\mathbb{R}}

\newcommand{\ZZ}{\mathbb{Z}}

\newtheorem{definition}{Definition}

\newtheorem{lemma}{Lemma}
\newtheorem{proposition}{Proposition}
\newtheorem{theorem}{Theorem}

\newtheorem{remark}{Remark}

\author{Emmanuel Briand}
\address{
Emmanuel Briand and Mercedes Rosas,
Departamento de \'Algebra,
Facultad de Matem\'aticas,
Aptdo. de Correos 1160,
41080 Sevilla, Spain.}
\email{ebriand@us.es, mrosas@us.es}

\author{Rosa Orellana}
\address{Rosa Orellana, Dartmouth College, Mathematics Department, 6188 Kemeny Hall, Hanover, NH 03755, USA.}
\email{rosa.c.orellana@dartmouth.edu}

\author{Mercedes Rosas}

\thanks{Emmanuel Briand is supported by a Juan de la Cierva Fellowship (MICINN, Spain). Mercedes Rosas is supported by a Ram\'on y Cajal Fellowship (MICINN, Spain). Both are also supported by Projects MTM2007--64509 (MICINN, Spain) and FQM333 (Junta de Andalucia).}

\title[Quasipolynomial formulas for Kronecker coefficients]{Quasipolynomial formulas for the Kronecker coefficients indexed by two two--row shapes\\ (extended abstract)}

\begin{document}

\begin{abstract}
We show that the Kronecker coefficients indexed by two two--row shapes are given by quadratic quasipolynomial formulas whose domains are the maximal cells of a fan. Simple calculations provide explicitly the quasipolynomial formulas and a description of the associated fan.

These new formulas are obtained from analogous formulas for the corresponding reduced Kronecker coefficients and a formula recovering the Kronecker coefficients from the reduced Kronecker coefficients.

As an application, we characterize all the Kronecker coefficients indexed by two two-row shapes that are equal to zero. This allowed us to  disprove a conjecture of Mulmuley about the behavior of the
 stretching functions attached to the Kronecker coefficients.

\medskip

\noindent {\bf R\'esum\'e.}
Nous démontrons que les coefficients de Kronecker indexés par deux partitions de longueur au plus 2 sont donnés par des formules quasipolynomiales quadratiques dont les domaines de validité sont les cellules maximales d'un éventail. 
Des calculs simples nous donnent une description explicite des formules quasipolynomiales et de l'éventail associé. 
Ces nouvelles formulas sont obtenues de formules analogues pour les coefficients de Kronecker réduits correspondants et au moyen d'une formule reconstruisant les coefficients de Kronecker à partir des coefficients de Kronecker réduits.

Une application est la caractérisation exacte de tous les coefficients de Kronecker non--nuls indexés par deux partitions de longueur au plus deux. Ceci nous a permis de réfuter une conjecture de Mulmuley au sujet des fonctions de dilatations associées aux coefficients de Kronecker.
\end{abstract}

\maketitle


\section*{Introduction}


A fundamental problem in algebraic combinatorics is the Clebsch-Gordan problem: given a linearly reductive group $G$, give a combinatorial description of the coefficients $m_{\mu\nu}^{\lambda}$ in the decomposition into irreducibles of the tensor product of two (finite-dimensional complex) irreducible representation $V_{\mu}(G)$ and $V_{\nu}(G)$:
\[V_\mu(G) \otimes V_\nu(G) \cong \bigoplus_{\lambda}m_{\mu\nu}^\lambda V_{\lambda}(G)\]
While this problem has been solved satisfactorily for the general linear group, $GL(n)$, the most elementary linear group, this is not the case for the symmetric group, $S_n$, the most fundamental finite group.  

In the case of $GL(n)$, the coefficients $m_{\mu\nu}^\lambda=c_{\mu\nu}^\lambda$ are the well known Littlewood-Richardson coefficients.  There exists several combinatorial descriptions for them.  One of these descriptions was given by \citet{Berenstein:Zelevinsky} that showed that $c_{\mu\nu}^\lambda$ counts the integral points in a well-defined family of polytopes.  This initiated a series of works concerning the stretching functions associated to these coefficients that culminated with the proof by \citet{Knutson:Tao:1} of the saturation conjecture. 
Finally,  \citet{Rassart}  showed that the Littlewood-Richardson coefficients $c_{\mu\nu}^{\lambda}$ are given by polynomial functions of the parts of $\lambda$, $\mu$ and $\nu$, on the maximal cells of a fan. 

For the symmetric group $S_n$, the coefficients $m_{\mu\nu}^\lambda = g_{\mu\nu}^\lambda$ are called the Kronecker coefficients.  Amazingly, there is no combinatorial description of these coefficients in general. Particular families have been investigated. In this paper the \emph{Kronecker coefficients indexed by two two--row shapes} are considered. They are the coefficients $g_{\mu\nu}^\lambda$ such that both $\mu$ and $\nu$ have two rows. Formulas for them have already been given by \citet{Remmel:Whitehead} and \citet{Rosas:2001}. 
Recent works by \citet{Luque:Thibon:4qubit,Garsia:Musiker:Wallach:Xin,Brown:VanWilli:Zabrocki} have revived the interest of obtaining better formulas for the Kronecker coefficients indexed by two two--row shapes as Hilbert series related to these coefficients have been linked to problems in quantum information theory.

New problems about the Clebsch--Gordan coefficients have been raised recently by the specialists of computational complexity. \citet{Narayanan} showed that the computation of the Littlewood--Richardson coefficients is a \#P--complete problem. \citet{Burgisser:Ikenmeyer} showed that the computation of the Kronecker coefficients is \#P--hard. On the other hand, the saturation property implies that the non--vanishing of a Littlewood--Richardson coefficient can be decided in polynomial time \citep{GCT3}. Is it also the case for the Kronecker coefficients? This question lies at the heart  of a detailed plan, \emph{Geometric Complexity Theory}, that \citet{GCT1:SIAM} elaborated to prove that $P \neq NP$ over the complex numbers (an arithmetic, non--uniform version of $P \neq NP$). This lead \citet{GCT6} to state a series of conjectures about the stretching functions associated to the Kronecker coefficients. The scarce information available about Kronecker coefficients made difficult even the experimental checking of these conjectures. By means of the formulas by \citet{Remmel:Whitehead} and \citet{Rosas:2001} it was only  possible to check them on large samples of Kronecker coefficients indexed by two two--row shapes \citep[see][]{GCT6}.

The present article obtains a new description for the Kronecker coefficients indexed by two two--row shapes, given by quasi--polynomial functions on the chambers of  fans, resembling the description of 
\citet{Rassart} for the Littlewood--Richardson coefficients. It is efficient enough to check Mulmuley's conjectures for all Kronecker coefficients indexed by two two--row shapes (and, actually, disprove them by providing explicit counter--examples). We start our investigation by looking at  Murnaghan's \emph{reduced Kronecker coefficients} $\bar{g}_{\alpha\beta}^\gamma$ \citep{Murnaghan:1938},  a related family of coefficients indexed by triples of partitions, which are stable values of stationary sequences of Kronecker coefficients.
Our first result expresses the Kronecker coefficients in terms of the reduced Kronecker coefficients (Theorem \ref{thm:gbar to g}). 
 Exploiting the work of \citet{Rosas:2001}  we are able to show that the reduced Kronecker polynomials  related to the two-row family count integral points in a polygon of $\mathbb{R}^2$.  From this we describe an explicit piecewise quasipolynomial formula for these reduced Kronecker coefficients. The pieces are the 26 maximal cells of a fan. Last, using our formula that recovers the Kronecker coefficients from the  reduced Kronecker coefficients, we obtain, with the help of the Maple package \emph{convex} by \citet{CONVEX}, explicit piecewise quasipolynomial formulas for the Kronecker coefficients indexed by two two-row shapes.  It is given by $74$ quadratic quasipolynomials whose domains are the maximal cells of a fan.   



As an application, we list all Kronecker coefficients indexed by two two-row shapes that are equal to zero. This made possible the discovery of counter--examples to Mulmuley's conjectures \citep{Briand:Orellana:Rosas:SH}.
In short, the advantage of our results is that for the first time we can completely study a complete nontrivial family of the Kronecker coefficients.

The detailed proofs will be presented in a full version \citep{Briand:Orellana:Rosas:Chamber} of this extended abstract.

\section{Piecewise Quasipolynomials}
\label{definicionquasi}

We now give a more detailed description of the main result.
A quasipolynomial is a function on $\ZZ^n$ given by polynomial formulas, whose domains are the cosets of a full rank sublattice of $\ZZ^n$. Remarkable examples of (univariate) quasipolynomials are the Ehrhart functions of polytopes of $\RR^k$ with rational vertices, that count the integral points in the dilations of the polytope \citep[see][chap. 4]{Stanley:vol1}.

We will obtain a description for the Kronecker coefficients indexed by two two--row shapes as a function of the following kind.
\begin{definition}\label{def}
A \emph{vector partition--like function} is a function $\phi$ on $\ZZ^n$ fulfilling the following: (i) There exists a convex rational polyhedral cone  $C$ such that $\phi$ is zero outside $C$. (ii) Inside $C$, the function $\phi$ is given by quasipolynomial formulas whose domains are (the sets of integral points of) the maximal (closed) cells of a fan $\fan$.

If $C$ and $\fan$ are as above and $Q$ is the family of quasipolynomial formulas, indexed by the maximal cells of $\fan$, we say that the triple $(C,\fan,Q)$ is \emph{a presentation of $\phi$ as a vector partition--like function}.
\end{definition}

\begin{remark}\label{rem:main}
A sum of vector partition--like functions $\phi_1$, $\phi_2$ is not necessarily vector partition--like. It is, however, the case when the functions admit presentations $(C,\fan,Q)$ and $(C',\fan',Q')$ with the same cone: $C=C'$. 
\end{remark}

Examples of vector partition--like functions are the vector partition functions, whose corresponding fans are the \emph{chamber complexes} \citep[see][]{Sturmfels, Brion:Vergne}.  

Vector partition--like functions also arise as functions counting integral solutions to some systems of linear inequalities depending on parameters.  Precisely, consider a system of inequalities of the form 
\begin{equation}\label{eq:system}
u_i(x)+c_i(h) \geq 0, \qquad i=1,\ldots, N
\end{equation}
 where the functions $u_i$ and $c_i$ are integral, homogeneous linear forms on $\RR^m$ and $\RR^n$ respectively. The unknown is $x$ and the parameter is $h$. 
 Assume that for any $h\in \RR^n$ the set of solutions $x$ of the system is bounded. Let $h \mapsto \phi(h)$ be the function that counts the \emph{integral} solutions $x$ of the system. This function $\phi$ is vector partition--like. This follows from the reduction of this function to a vector partition function \citep[see][]{Brion:Vergne}. Here the cone $C$ in Definition \ref{def} is the set of values of the parameter $h$ making the system feasible.

Let $\ell$ be a positive integer. The function $(\lambda,\mu,\nu) \mapsto c_{\mu,\nu}^{\lambda}$  from triples of partitions with at most $\ell$ parts to Littlewood--Richardson coefficients is vector partition--like. This is because this function counts the integral solutions of a system of inequalities depending on parameters (the parts of the partitions) of the form \eqref{eq:system}. Indeed, such a system can be derived from the Littlewood--Richardson rule \cite[see][]{GCT3}. Alternatively, one can use the system defining Knutson and Tao's Hive polytopes \citep[see the exposition by][]{Buch}. 

It is natural to ask if similar results also hold for the Kronecker coefficients. Let $\ell_1$ and $\ell_2$ be positive integers. If $\mu$ and $\nu$ are partitions of length at most $\ell_1$ and $\ell_2$ respectively then $g_{\mu,\nu}^{\lambda}$ can be nonzero only if $\lambda$ has at most $\ell_1 \ell_2$ parts. 
The analogous function to consider is thus $G_{\ell_1,\ell_2} : (\lambda,\mu,\nu)\mapsto g_{\mu,\nu}^{\lambda}$ defined on triples of partitions with at most $\ell_1 \ell_2$, $\ell_1$ and $\ell_2$ parts respectively.  No interpretation of the functions $G_{\ell_1,\ell_2}$ as counting integral solutions to systems of inequalities of the form \eqref{eq:system} is known.  Nevertheless, very close results were obtained by \cite{GCT6}: (i) The functions $G_{\ell_1,\ell_2}$ fulfill the conditions in Definition \ref{def} with $\fan$ a complex of polyhedral cones instead of a fan. (ii) For any $\lambda$, $\mu$, $\nu$, the \emph{stretching function} $N \in \NN \mapsto g_{N\mu,N\nu}^{N^{\lambda}}$ is a univariate quasipolynomial. Here $N\lambda$ stands for the partitions obtained from $\lambda$ by multiplying all parts by $N$. Combining these two results, one gets that the functions $G_{\ell_1,\ell_2}$ fulfill the conditions in the definition of vector partition--like with ``maximal closed cells'' replaced with ``open cells'' in (ii).

The simplest non--trivial case is $G_{2,2}$, describing the Kronecker coefficients indexed by two two--row shapes. Even this case is somehow difficult. In this work we prove the following:
\begin{theorem}\label{thm:main1}
The function 
\[
G_{2,2}:\;(\lambda_1,\ldots,\lambda_4,\mu_1,\mu_2,\nu_1,\nu_2)\in \ZZ^8 \mapsto g_{(\mu_1,\mu_2)(\nu_1,\nu_2)}^{(\lambda_1,\lambda_2,\lambda_3,\lambda_4)}
\]
is vector partition--like.
\end{theorem}
\begin{remark}\label{rem:G}
A Kronecker coefficient $g_{\mu,\nu}^{\lambda}$ can be nonzero only if its three indexing partitions have the same weight. This and the formula
$g_{(\mu_1,\mu_2)(\nu_1,\nu_2)}^{(\lambda_1,\lambda_2,\lambda_3,\lambda_4)}=g_{(\mu_1-2,\mu_2-2)(\nu_1-2,\nu_2-2)}^{(\lambda_1-1,\lambda_2-1,\lambda_3-1,\lambda_4-1)}$ reduce the study of $G_{2,2}$ to the study of the function 
\[
(n,\gamma_1,\gamma_2,r,s) \mapsto g_{(n-r,r)(n-s,s)}^{(n-\gamma_1-\gamma_2,\gamma_1,\gamma_2)}
\]
\end{remark}

\section{Murnaghan's Theorem and reduced Kronecker coefficients}\label{Murnagham}

In this section we introduce Murnaghan's reduced Kronecker coefficients $\overline{g}_{\alpha,\beta}^{\gamma}$. They are integers indexed by triples of partitions closely related to the Kronecker coefficients. The Kronecker coefficients indexed by two two--row shapes will be re--obtained from the reduced Kronecker coefficients indexed by two one--row shapes (Section \ref{sec:from gbar to g}) which will be easy to describe (Theorem \ref{cor:Pi} and Section \ref{sec:formulas red}).

The Jacobi--Trudi formula expresses the Schur functions as determinants in the complete sums $h_k$. When $\lambda$ has at most $k$ parts, it asserts that:
\[
{s}_{\lambda}=\det( h_{j-i+\lambda_i})_{i,j=1,\ldots,k}
\]
(where $h_k=0$ when $k<0$, $h_0=1$ and $\lambda_i=0$ for $i$ greater than the length of $\lambda$.)

This formula can also be applied in the case when $\lambda$ is not a partition, i.e. is not nondecreasing. The functions $s_{\lambda}$ obtained are either $0$, or Schur functions up to a sign.

Let $n$ be an integer and $\lambda$ a partition. Then $|\lambda|$ stands for the sum of the parts of $\lambda$ and for any integer $n$, we denote with $(n-|\lambda|,\lambda)$ the sequence $(n-|\lambda|,\lambda_1,\lambda_2,\ldots)$. This is a partition if and only if $n \geq |\lambda|+\lambda_1$. Last $\cut{\lambda}$ stands for the partition $(\lambda_2,\lambda_3,\ldots)$, which is obtained by removing the first part of $\lambda$.

\begin{theorem}[\citet{Murnaghan:1938,Murnaghan:1955}]\label{thm:stability}
There exists a family of nonnegative integers $(\overline{g}_{\alpha,\beta}^{\gamma})$ indexed by triples of partitions $(\alpha,\beta,\gamma)$ such that, for fixed partitions $\alpha$ and $\beta$, only finitely many terms $\overline{g}_{\alpha,\beta}^{\gamma}$ are non--zero, and 
for all $n \geq 0$, 
\[
{s}_{(n-|\alpha|,\alpha)} \ast {s}_{(n-|\beta|,\beta)}
=
\sum_{\gamma} \overline{g}_{\alpha,\beta}^{\gamma} {s}_{(n-|\gamma|,\gamma)}
\]
\end{theorem}
Following \citet{Klyachko}, we call the coefficients $\overline{g}_{\alpha,\beta}^{\gamma}$ the \emph{reduced Kronecker coefficients}. They are called \emph{extended Littlewood--Richardson numbers} in \citet{Kirillov:saturation}  because of the following property, observed first in \citet{Murnaghan:1955} and proved in \citet{Littlewood:1958}: if $\alpha$, $\beta$ and $\gamma$ are three partitions such that $|\gamma|=|\alpha|+|\beta|$ then $\overline{g}_{\alpha\beta}^{\gamma}=c_{\alpha\beta}^{\gamma}$.

\begin{remark}
It follows from Murnaghan's Theorem that for fixed partitions $\alpha$, $\beta$, $\gamma$, the sequence of Kronecker coefficients $g_{(n-|\alpha|,\alpha),(n-|\beta|,\beta)}^{(n-|\gamma|,\gamma)}$ ($n$ big enough so that all three indices are partitions) is stationary with limit $\overline{g}_{\alpha\beta}^{\gamma}$.
\end{remark}

\section{From reduced to non--reduced Kronecker coefficients}\label{sec:from gbar to g}

In this section we give a formula that allows us to recover the non-reduced Kronecker coefficients from the reduced Kronecker coefficients, and we apply it for the Kronecker coefficients indexed by two two--row shapes.

\begin{theorem}\label{thm:gbar to g}
Let $\ell_1$, $\ell_2$ and $n$ be positive integers. Let 
$\lambda$, $\mu$, $\nu$ be partitions of $n$ such that $\mu$ has length at most $\ell_1$ and $\nu$ has length at most $\ell_2$.
If $\lambda$ has length at most $\ell_1 \ell_2$ then
\begin{equation}\label{eq:g to gbar}
g_{\mu \nu}^{\lambda}=
\sum_{i=1}^{\ell_1 \ell_2} (-1)^{i+1} \overline{g}_{\cut{\mu}, \cut{\nu}}^{\lambda^{\dagger i}}
\end{equation}
where $\lambda^{\dagger i}$ is the partition obtained from $\lambda$ by incrementing the $i-1$ first parts and removing the $i$--th part, that is:
\[
\lambda^{\dagger i}=\left(1+\lambda_1,1+\lambda_2,\ldots,1+\lambda_{i-1},\lambda_{i+1},\lambda_{i+2}, \ldots\right)
\]
\end{theorem}

Formula \eqref{eq:g to gbar} applies as follows in the case $\ell_1=\ell_2=2$ and $\lambda$ with at most three parts:
\begin{equation}\label{eq:g to gbar22}
g_{(n-r,r)(n-s,s)}^{(\lambda_1,\lambda_2,\lambda_3)} = 
\overline{g}_{(r)(s)}^{(\lambda_2,\lambda_3)}
-
\overline{g}_{(r)(s)}^{(\lambda_1+1,\lambda_3)}
+
\overline{g}_{(r)(s)}^{(\lambda_1+1,\lambda_2+1)}
\end{equation}
where $n=|\lambda|$.
(One can show that the last expected term, $\overline{g}_{(r)(s)}^{(\lambda_1+1,\lambda_2+1,\lambda_3+1)}$, is always zero.) 

The reduced Kronecker coefficients that appear in this formula are all of the form $\overline{g}_{(r)(s)}^{(\gamma_1,\gamma_2)}$. These coefficients admit the following description, derived in \citet{Briand:Orellana:Rosas:SH} from the description for the Kronecker coefficients indexed by two two--row shapes provided by \citet{Rosas:2001}. An equivalent description for the reduced Kronecker coefficients indexed by two one--row shapes is given by \citet{Thibon}.
\begin{theorem}[\citet{Briand:Orellana:Rosas:SH}]\label{cor:Pi}
Let $r$, $s$ and $\gamma_1 \geq \gamma_2$ be nonnegative integers and $h=(r,s,\gamma_1,\gamma_2)$. The reduced Kronecker coefficient $\overline{g}_{(r)(s)}^{(\gamma_1,\gamma_2)}$ counts the integral solutions to the system of inequalities $u_i(X,Y) + c_i(h) \geq 0$ for $i=0$, \ldots, $6$, where:
\begin{equation}\label{eq:u}
\begin{array}{l@{\quad}l}
\begin{array}{rcl}
u_0(v) + c_0(h)&=& X-s\\
u_1(v)+c_1(h)&=& X-r\\
u_2(v)+ c_2(h)&=& X+Y-r-s+\gamma_1 \\
u_3(v)+ c_3(h)&=& Y 
\end{array}
&
\begin{array}{rcl}
u_4(v)+ c_4(h)&=& Y-X+|\gamma| \\
u_5(v)+ c_5(h)&=& -X-Y+r+s-\gamma_2\\
u_6(v)+ c_6(h)&=& X-Y-\gamma_1
\end{array}
\end{array}
\end{equation}
In particular, the function $R: (r,s,\gamma_1,\gamma_2) \in \ZZ^4 \mapsto \overline{g}_{(r)(s)}^{(\gamma_1,\gamma_2)}$ is vector partition--like.
\end{theorem}

Theorem \ref{cor:Pi} and Formula \eqref{eq:g to gbar22} provide a piecewise quasipolynomial description for $G_{2,2}$ (see Remark \ref{rem:G}). But the corresponding domains of quasipolynomiality obtained are neither closed, nor cones. The remainder of this work is devoted to correct this and obtain, still from Theorem \ref{cor:Pi} and Formula \eqref{eq:g to gbar22} a vector partition--like presentation for $G_{2,2}$.

The main tools are the Lemma \ref{lemma:Tv}, below, and an explicit vector partition--like presentation for the function $R$ (section \ref{sec:formulas red}) showing that the lemma applies.

Let $\FF_0$, $\FF_1$, $\FF_2$ be the linear maps from $\RR^5$ to $\RR^4$ that send $(n,r,s,\gamma_1,\gamma_2)$ to $(r,s,\gamma_1,\gamma_2)$, $(r,s,n-\gamma_1-\gamma_2,\gamma_2)$, $(r,s,n-\gamma_1-\gamma_2,\gamma_1)$ respectively. Let $T_1$ and $T_2$ be the translations in $\RR^4$ of vector $v_1=(0,0,1,0)$ and $v_2=(0,0,1,1)$ respectively. 

Let $\SK$ (resp. $\RK$) be the cone of $\RR^5$ (resp. of $\RR^4$) generated by all $(n,r,s,\gamma_1,\gamma_2) \in \ZZ^5$ (resp. all $(r,s,\gamma_1,\gamma_2) \in \ZZ^4$) such that  the Kronecker coefficient $g_{(n-r,r)(n-s,s)}^{(n-\gamma_1-\gamma_2,\gamma_1,\gamma_2)}$ (resp. the reduced Kronecker coefficient $\overline{g}_{(r)(s)}^{(\gamma_1,\gamma_2)}$) is defined and positive. The explicit description of $\SK$ is provided by \citet{Bravyi} (see also the general approach by \citet{Klyachko}). The cone $\RK$ is the image of $\SK$ under $F_0$.

For $x \in \ZZ^5$ set $\chi_{\SK}(x)=1$ if $x \in \SK$ and $\chi_{\SK}(x)=0$ otherwise. Then we can rewrite Formula \ref{eq:g to gbar22} as follows:
\[
G(x)=R \circ \FF_0(x) - \chi_{\SK}(x) \cdot R \circ T_1 \circ \FF_1(x)
 + \chi_{\SK}(x) \cdot R \circ T_2 \circ \FF_2(x)
\]
where $G(x)=G(n,r,s,\gamma_1,\gamma_2)=g_{(n-r,r)(n-s,s)}^{(n-\gamma_1-\gamma_2,\gamma_1,\gamma_2)}$ when $(n-r,r)$, $(n-s,s)$, $(n-\gamma_1-\gamma_2,\gamma_1,\gamma_2)$ are partitions, and $G(n,r,s,\gamma_1,\gamma_2)=0$ otherwise.

After Remark \ref{rem:main}, Theorem \ref{thm:main1} will be proved if we show that all three vector partition--like functions $R \circ \FF_0$, $\chi_{\SK} \cdot R \circ T_1 \circ \FF_1$ and $\chi_{\SK} \cdot R \circ T_2 \circ \FF_2 $ admit presentations with the same cone: $(\SK,\fan_0,Q_0)$, $(\SK,\fan_1,Q_1)$ and $(\SK,\fan_2,Q_2)$. 

That $R \circ \FF_0$ admits a presentation $(\SK,\fan_0,Q_0)$ is immediate because $F_0^{-1}(\RK)=\SK$. To show that $\chi_{\SK} \cdot R \circ T_1 \circ \FF_1$ and $\chi_{\SK} \cdot R \circ T_2 \circ \FF_2 $ also admit presentations with cone $\SK$ we will need to apply two times Lemma \ref{lemma:Tv} below, with $p=5$, $q=4$, $C=\SK$, $C'=\RK$, $\phi=R$ and $F=F_i$, $v=v_i$ for $i=1$, $2$.

Given subsets $A$, $B$ of $\RR^q$ we denote with $A+B$ the set $\{a+b\,|\,a \in A, b \in B\}$. Given $v \in \RR^q$ and $I$ subset of $\RR$ we denote with $I\,v$ the set $\{x v\,|\, x \in I\}$.
\begin{lemma}\label{lemma:Tv}
Let $\phi$ be a vector partition--like function on $\ZZ^q$ with presentation $(C',\fan',Q)$. Let $C$ be a convex rational polyhedral cone of $\RR^p$ and $F$ an integral linear map  from $\RR^p$ onto $\RR^q$. Let $v \in \ZZ^q$ and $T$ be the translation of $\RR^q$ of vector $v$.
Let $\fan$ be the fan subdividing $C \cap F^{-1}(C')$, whose cells are all sets of the form $C \cap F^{-1}(\sigma')$ for $\sigma'$ cell of $\fan'$.

Assume that the cone $C \cap F^{-1}(C')$ is full--dimensional in $\RR^p$. Assume also that:
\begin{itemize}
\item[(a)] Whenever $H$ is a hyperplane separating two adjacent maximal cells $\sigma'_1$, $\sigma'_2$ of $\fan'$ such that $F(C)$ is not included in $H+\RR_+ v$, the following holds: The quasipolynomials $Q_{\sigma'_1}$ and $Q_{\sigma'_2}$ coincide on the integral points of the strip $H\,+\,]0;1]\,v$.
\item[(b)] Whenever $H$ is a hyperplane containing a facet of $C'$, such that $\RR_+\,v\,+\,F(C)$ is not contained in the half--plane $H+C'$, the following holds: For all maximal cells $\sigma'$ of $\fan'$ having a facet contained in $H$, the quasipolynomial $Q_{\sigma'}$ vanishes on the integral points of the strip $H\,+\,]0;1]\,v$.
\end{itemize}
Then
\begin{itemize}
\item[(i)] The function $\phi \circ T \circ F$ is zero on the integral points of the closure of $C \setminus F^{-1}(C')$.
\item[(ii)] If $C \cap F^{-1}(\sigma')$ is a maximal cell of $\fan$ (where $\sigma'$ is a maximal cell of $\fan'$) then $\phi \circ T \circ F$ and $Q_{\sigma} \circ T \circ F$ coincide on its integral points.
\end{itemize}
\end{lemma}

Applying the lemma as indicated requires a precise description of a presentation $(\RK,\fan_R,Q_R)$ of $R$.  The next section provides such a description.

\section{Formulas for the reduced Kronecker coefficients indexed by two one--row shapes}\label{sec:formulas red}

Let $u_i$ and $c_i$, for $i=0$, $1$, \ldots, $6$ be the integral linear forms defined in \eqref{eq:u}. After \citet{Brion:Vergne}, the function $\psi$ that associates to $y \in \ZZ^7$ the number of integral solutions  of the system $u_i(X,Y)+y_i \geq 0$, $i=0$, \ldots, $6$ is a vector partition function. In particular, it admits a very well--described vector partition--like presentation $(C_{\psi},\mathcal{F}_{\psi},Q_{\psi})$. The corresponding fan is the \emph{chamber complex} of $\psi$, see \citet{Brion:Vergne,Sturmfels}. 

Remember (Theorem \ref{cor:Pi}) that $R$ is the function that associates the reduced Kronecker coefficient $\overline{g}_{(r)(s)}^{(\gamma_1,\gamma_2)}$ to $(r,s,\gamma_1,\gamma_2) \in \ZZ^4$. Then $R=\psi \circ c$, where $c$ is the linear map from $\RR^4$ to $\RR^7$ that maps $h=(r,s,\gamma_1,\gamma_2)$ to $(c_0(h),c_1(h),\ldots,c_6(h))$. Therefore, one obtains a very explicit vector partition--like presentation $(c^{-1}(C_{\psi}),\mathcal{\fan}_R,Q_R)$ for $R$ by taking for ${\fan}_R$ the inverse image of ${\fan}_{\psi}$ under $c$, and for $Q_R$ the family of functions $Q_{R,c^{-1}(\sigma)}=Q_{\psi,\sigma}\circ c $ for $\sigma$ maximal cell of $\fan_{\psi}$. We present this description.

Let $h \in \RR^7$. Denote with $\Pi(h)$ the set of real solutions of the system \eqref{eq:u}. For $i=0$, $1$, \ldots, $6$, let $L_i(h)$ be the line with equation $a_i X + b_i Y +c_i(h)=0$ where $u_i(X,Y)=a_i X + b_i Y$. 

For any three elements $i$, $j$, $k$  of $\{0,1,\ldots,6\}$ define:
\begin{equation}\label{eq:fijk}
f_{ijk}(h)=
-\left|
\begin{matrix}
a_i & a_j & a_k \\
b_i & b_j & b_k \\
c_i(h) & c_j(h) & c_k(h) 
\end{matrix}
\right|
\end{equation}
Define also $f_{25}=\gamma_1-\gamma_2$ and $f_{46}=\gamma_2$. The linear form  $f_{25}$ (resp. $f_{46}$) is proportional to $f_{25k}$ for all $k \neq 2,5$ (resp.: to $f_{46k}$ for all $k \neq 4,6$) and 
its vanishing is the condition for the two parallel lines $L_2$ and $L_5$ (resp. $L_4$ and $L_6$) to coincide.

\begin{figure}[htbp]
\[
\xymatrix{
& 
\text{\small 1245}\ar@{-}[rr]            
&& 
\text{\small 12456}\ar@{-}[rrrr]
&&&
&\text{\small 02456}
&&
\text{\small 0245}\ar@{-}[ll] 
\\
\text{\small 145}\ar@{-}[rr]\ar@{-}[ru]\ar@{-}[dd]
&& 
\text{\small 1456}\ar@{-}[ru]\ar@{-}[dd]
&&&&
\text{\small 0456}\ar@{-}[ru]\ar@{-}[dd]
&&
\text{\small 045}\ar@{-}[ll]\ar@{-}[ru]\ar@{-}[dd]
&\\&
\text{\small 12345}\ar@{-}[uu]|
\hole
\ar@{-}[rr]|!{[ur];[dr]}\hole
&& 
\text{\small 123456}\ar@{-}[uu]\ar@{-}[rrrr]|!{[drrr];[urrr]}\hole
&&&&
\text{\small 023456}\ar@{-}[uu]|
\hole
&&
\text{\small 02345}\ar@{-}[uu]\ar@{-}[ll]|!{[ul];[dl]}\hole 
\\
\text{\small 1345}\ar@{-}[rr]\ar@{-}[ru]             
&&
\text{\small 13456}\ar@{-}[ru]               
&&
\text{\small 3456}\ar@{-}[ll]\ar@{-}[rr]
&&
\text{\small 03456}\ar@{-}[ru]
&&
\text{\small 0345}\ar@{-}[ll]\ar@{-}[ru]
&\\&
\text{\small 1235}\ar@{-}[uu]|
\hole\ar@{-}[rr] |!{[dr];[ur]}\hole     
&&
\text{\small 12356} \ar@{-}[uu] |
\hole\ar@{-}[rrrr]|!{[dr];[ur]}\hole|!{[drrr];[urrr]}\hole
&&&&
\text{\small 02356}\ar@{-}[uu]|
\hole
&&
\text{\small 0235}\ar@{-}[uu]\ar@{-}[ll]|!{[dl];[ul]}\hole
\\
\text{\small 135}\ar@{-}[ru]\ar@{-}[uu]\ar@{-}[rr]
&& 
\text{\small 1356}\ar@{-}[ru]\ar@{-}[uu] 
&&
\text{\small 356}\ar@{-}[ll]\ar@{-}[uu]\ar@{-}[rr]
&& 
\text{\small 0356}\ar@{-}[ru]\ar@{-}[uu] 
&& 
\text{\small 035}\ar@{-}[ru]\ar@{-}[uu]\ar@{-}[ll]
&
}
\]
\caption{The graph $\mathcal{G}$.}\label{fig:graph}
\end{figure}
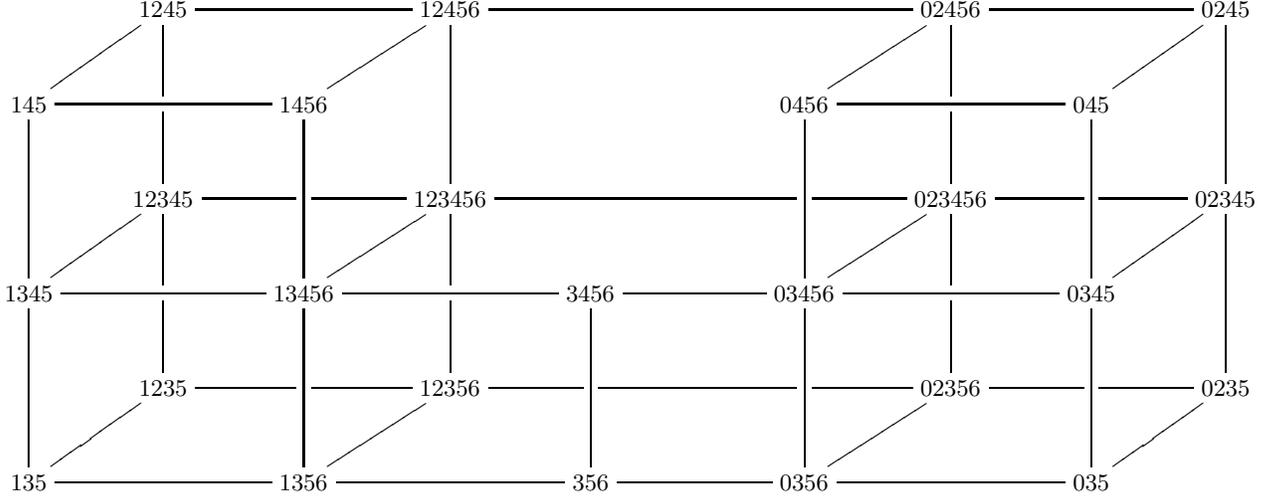
\begin{itemize}
\item {\bf The cone $c^{-1}(C_{\psi})$} is equal to the cone $\RK$ introduced in Section \ref{sec:from gbar to g}. It is defined by the system of linear inequalities:
\[
f_{145} \leq 0,\quad
f_{045} \leq 0,\quad
f_{356} \leq 0,\quad
f_{035} \leq 0,\quad
f_{135} \leq 0,\quad
f_{25} \geq 0,\quad
f_{46} \geq 0.
\]
\item {\bf The fan $\fan_R$}: Let $S$ be the locus of parameters $h$ such that three lines $L_i(h)$, $L_j(h)$, $L_k(h)$ meet in $\Pi(h)$. 
The fan $\fan_R$ is the fan whose chambers (maximal open cells) are the connected components of $\RK \setminus S$. In each chamber $\sigma$ the set of indices $i$ such that $L_i(h)$ supports a side of $\Pi(h)$ is constant. Denote this set with $\textsf{Sides}(\sigma)$. This set $\textsf{Sides}(\sigma)$ determines $\sigma$. Therefore we denote a chamber $\sigma$ with $\sigma_I$ when $\textsf{Sides}(\sigma)=I$, e.g. $\sigma_{1245}$ for the chamber $\sigma$ such that $\textsf{Sides}(\sigma)=\{1,2,4,5\}$. There are $26$ chambers $\sigma_I$ in $\fan_R$. The corresponding indices $I=\textsf{Sides}(\sigma_I)$ are the vertices of the graph $\mathcal{G}$ in Figure \ref{fig:graph}. Adjacency in $\mathcal{G}$ represents adjacency in $\fan_R$: chambers $\sigma_I$ and $\sigma_J$ are adjacent (i.e. their closures have a common facet) if and only if $I$ and $J$ are adjacent vertices in $\mathcal{G}$. Observe that when $\sigma_I$ and $\sigma_J$ are adjacent then:
\begin{itemize}
\item either $I$ and $J$ are obtained from each other by exchanging $0$ and $1$. Then $\sigma_I$ and $\sigma_J$ are separated by the hyperplane of equation $r=s$. There is $r>s$ on $\sigma_I$ if $1 \in I$.
 \item or one of the sets is obtained from the other by inserting a unique element. Say $J=I \cup \{j\}$ with $j \not \in I$. If the elements of $J$ are $p_1<p_2<\cdots<p_t$ say that the successor of $p_q$ is $p_{q+1}$, for $q=1$, \ldots, $t-1$, and that the successor of $p_t$ is $p_1$. This defines a cyclic order on $J$. Let $i$ and $k$ be the predecessor and successor of $j$ in this cyclic order. Then $\sigma_I$ and $\sigma_J$ are separated by the hyperplane of equation $f_{ijk}=0$, and $f_{ijk}>0$ on $\sigma_I$. 
\end{itemize} 
\item {\bf The quasipolynomial formulas on each maximal cell:} For simplicity we set $q_I=Q_{R,\overline{\sigma_I}}$. This is the quasipolynomial formula for $R$ valid on the cell $\overline{\sigma_I}$ (the topological closure of the chamber $\sigma_I$). Rather than displaying explicit expressions for all quasi--polynomials $q_I$, it is enough to present one of them (we choose $q_{135}$) and display all differences $q_I-q_J$ for $\sigma_I$ and $\sigma_J$ adjacent. All quasi--polynomials $q_I$ can be recovered easily from this information by chasing on the graph $\mathcal{G}$ (Figure \ref{fig:graph}), \emph{e.g.} 
\[
q_{1456}=(q_{1456}-q_{1{\bf 3}456})+(q_{13{\bf 4}56}-q_{1356})+(q_{135{\bf 6}}-q_{135})+q_{135}
\] 
There is:
\[
q_{135}(r,s,\gamma_1,\gamma_2)=
\frac{1}{2}\,
\left(
s-\gamma_2+1 
\right)
\left(
s-\gamma_2+2
\right)
\]
Let $\sigma_I$ and $\sigma_J$ be two adjacent chambers of $\fan$.
\begin{itemize}
\item If $I$ and $J$ are obtained from each other by exchanging $0$ and $1$ then $q_I=q_J$. 
\item If $J=I \cup \{j\}$ with $j \not \in I$ then $q_I-q_J$ depend only of $j$ and its predecessor $i$ and successor $k$ in $J$, and is as indicated in Table \ref{table:differences}.
\end{itemize}
\end{itemize}
\begin{table}
\[
\begin{array}{|c|c|c|}
\hline
ijk & q_I(h)-q_J(h) & 
\begin{array}{c}
\text{\bf Values $\delta$ such that}\\
\text{\bf $q_I$ and $q_J$ coincide}\\
\text{\bf on $f_{ijk}=\delta$}
\end{array}
\\\hline\hline
\begin{array}{cc}613, 123, 134 \\ 603, 023, 034\end{array}& \frac{1}{2} f_{ijk}(h) \left(f_{ijk}(h)-1\right) & 0,1 \\
\hline
234 & \frac{1}{4} \left(f_{ijk}(h)\right)^2 + 
\left \lbrace
\begin{array}{cl}
0 &\text{ if } f_{ijk}(h) \equiv 0 \mod 2 \\
-1/4 &\text{ else.}
\end{array}
\right.& -1,0,1\\
\hline
\begin{array}{c}345, 124, 561 \\ 024, 560 \end{array} &
\frac{1}{4} f_{ijk}(h) \left(f_{ijk}(h)-2\right) + 
\left \lbrace
\begin{array}{cl}
0   &\text{ if } f_{ijk}(h) \equiv 0 \mod 2 \\
1/4 &\text{ else}
\end{array}
\right.
& 0,1,2 \\\hline
\end{array}
\]
\caption{The differences $q_I-q_J$ for $\sigma_I$ and $\sigma_J$ adjacent chambers of $\fan$.}\label{table:differences}
\end{table} 


If $\sigma_I$ and $\sigma_J$ are adjacent, the quasi--polynomials $q_I$ and $q_J$ coincide not only on the affine hyperplane spanned by the facet $\overline{\sigma_I} \cap \overline{\sigma_J}$ but also on close parallel hyperplanes.

\begin{proposition}\label{cor:coincidence interne}
Let $\sigma_I$ and $\sigma_J$ be two adjacent chambers of $\fan$ such that $J =I \cup \{j\}$ with $j \not \in I$. Let $i$ and $k$ be the predecessor and successor, respectively, of $j$ in $J$. 

Then $q_I-q_J$ coincide on the affine hyperplanes $f_{ijk}=\delta$ for the values of $\delta$ given by the third column in Table \ref{table:differences}.
\end{proposition}

Similarly, if the hyperplane $H$ supports a facet of a maximal cell $\overline{\sigma_I}$, and this facet is contained in the border of $\RK$, then $q_I$ vanishes on affine hyperplanes close and parallel to $H$.

\begin{proposition}\label{cor:coincidence externe}
Let $\sigma_I$ be a chamber of $\fan$ and $\tau$ an external facet of $\overline{\sigma_I}$ (i.e. a facet contained in the border of $\RK$). The hyperplane supporting $\tau$ admits as equation $f=0$ where $f$ is one of the linear forms $f_{145}$, $f_{045}$, $f_{356}$, $f_{035}$, $f_{135}$, $f_{25}$, $f_{46}$.

The set of values $\delta \in \ZZ$ such that $f$ vanishes identically on the affine hyperplane of equation $f=\delta$ is provided by  Table \ref{table:external}.
\end{proposition}

\begin{table}
\[
\begin{array}{|c|c|c|}
\hline
 \text{\bf Form $f$} 
& \multicolumn{1}{|p{5cm}|}{\text{\bf Chambers having a facet}

 \text{\bf supported by $f=0$}}
&
\begin{array}{c}
\text{\bf Values $\delta$ such that}\\
\text{\bf $q_I$ vanishes identically}\\
\text{\bf on $f=\delta$}
\end{array}
\\
\hline\hline
f_{46}= \gamma_2& 
{3456},{1456}, {0456} 
&
-1
\\ \hline
f_{25}= \gamma_1-\gamma_2 & 
{1245}, {0245},{1235}, {0235}
&
-1 
\\ \hline
f_{145}= r-s-\gamma_1& 145  & 1,2,3
\\ \hline
f_{045}= s-r-\gamma_1& {045} & 1,2,3
\\ \hline
f_{356}= |\gamma|-r-s & {356} & 1,2,3
\\ \hline
f_{035}= \gamma_2-r& {035}  & 1,2
\\ \hline
f_{135}= \gamma_2-s& {135}  & 1,2
\\ \hline
\end{array}
\]
\caption{The linear forms defining the facets of $\RK$.}\label{table:external}
\end{table}

It is immediate that $R \circ F_0$ has a vector partition--like presentations $(\SK,\fan_0,Q_0)$. Propositions \ref{cor:coincidence interne} and \ref{cor:coincidence externe} are used to apply Lemma \ref{lemma:Tv} and show that $\chi_{\SK} \cdot R \circ T_1 \circ \FF_1$ and $\chi_{\SK} \cdot R \circ T_2 \circ \FF_2$ have vector partition--like presentations $(\SK,\fan_1,Q_1)$ and $(\SK,\fan_2,Q_2)$. After Remark \ref{rem:main}, this proves Theorem \ref{thm:main1} and provides a way to compute a vector partition--like presentation for $G$ and  $G_{2,2}$.

\section{Formulas for the Kronecker coefficients indexed by two two--row shapes}\label{sec:Applications}

Once the presentations $(\SK,\fan_0,Q_0)$, $(\SK,\fan_1,Q_1)$, $(\SK,\fan_2,Q_2)$ for $R \circ F_0$, $\chi_{\SK} \cdot R \circ T_1 \circ F_1$ and $\chi_{\SK} \cdot R \circ T_2 \circ F_2$ have been determined, an explicit presentation $(\SK,\fan_{3},Q_3)$ for $G$ is obtained: The cells of $\fan_{3}$ are the intersection $\sigma_0 \cap \sigma_1 \cap \sigma_2$ for $\sigma_i$ a cell of $\fan_i$, $i\in\{0,1,2\}$. If $\sigma_0 \cap \sigma_1 \cap \sigma_2$ is a maximal cell of $\fan_3$ then the corresponding quasipolynomial formula for $G$ is $Q_{0,\sigma_0}-Q_{1,\sigma_1}+Q_{2,\sigma_2}$. We computed the description for $\fan_3$ by using the Maple Package CONVEX by \citet{CONVEX}: it has 177 maximal cells. It turns out that on some of them $G$ is given by the same quasipolynomial formulas, and that they can be glued together to form the maximal cells of a new fan $\fan_K$. In the new presentation $(\SK,\fan_K,P)$ obtained for $G$ the fan $\fan_K$ has only $74$ maximal cells. 

All $74$ quasipolynomial formulas $P_{\sigma}$ have the following form:
\begin{equation}\label{eq:shape}
P_{\sigma}
=1/4\;Q_{\sigma}+1/2\;L_{\sigma}+M_{\sigma}/4
\end{equation}
 where $Q_{\sigma}$ and $L_{\sigma}$ are
integral homogeneous polynomials in $(n,r,s,\gamma_1,\gamma_2)$
 respectively quadratic and linear. The function $M_{\sigma}$ takes integral values, fulfills $M_{\sigma}(0)/4=1$ and is constant on each coset of $\ZZ^5$ modulo the sublattice defined by $r+s \equiv n \equiv \gamma_1 \equiv \gamma_2 \equiv 0 \mod 2$.

Moreover, for all maximal cells $\sigma$, the functions  $Q_{\sigma}$, $L_{\sigma}$ are nonnegative on $\sigma$. This also holds for $M_{\sigma}$, for all cells $\sigma$ except four. This makes specially easy studying the \emph{support} of the Kronecker coefficients indexed by two two--row shapes. This is the set of all triples $(\lambda,\mu,\nu)$ such that $g_{\mu,\nu}^{\lambda}>0$ and $\mu$ and $\nu$ have at most two parts.

We obtain the following result. Let $(n,r,s,\gamma_1,\gamma_2) \in \SK$. Then $g_{(n-r,r)(n-s,s)}^{(n-\gamma_1-\gamma_2,\gamma_1,\gamma_2)}$ is zero if and only if at least one of the following five systems of conditions is fulfilled:
\begin{equation}\label{eq:conditions}
\begin{array}{lll}
\begin{array}{l}
\left\lbrace
\begin{array}{l}
n=2\,s=2\,r\\
\gamma_1 \text{ or } \gamma_2 \text{ odd.}
\end{array}
\right.
\\
\\
\left\lbrace
\begin{array}{l}
n=\max(2\,r,2\,s)\\
\gamma_1=\gamma_2\\
r+s+\gamma_1 \text{ odd.}
\end{array}
\right.
\end{array}
&
\begin{array}{l}
\left\lbrace
\begin{array}{l}
n=\max(2\,r,2\,s, |\gamma|+\gamma_1)\\
\gamma_2=0\\
r+s+\gamma_1 \text{ odd.}
\end{array}
\right.
\\
\\
\left\lbrace
\begin{array}{l}
n=|\gamma|+\gamma_1=\max(2\,r,2\,s)\\
r+s+\gamma_1 \text{ odd.}
\end{array}
\right.
\end{array}
&
\left\lbrace
\begin{array}{l}
n=max(2\,r,2\,s)\\
|r-s|=1\\
\min(2\,r,2\,s) \geq |\gamma|+\gamma_1\\
\gamma_1 \text{ or }\gamma_2 \text{ even.}
\end{array}
\right.
\end{array}
\end{equation}
This exhaustive description led us to a family of counterexamples for SH, a saturation conjecture formulated by \citet{GCT6}. The \emph{stretching functions} $\widetilde{g}_{\mu,\nu}^{\lambda}: N \mapsto g_{N\mu,N\nu}^{N\lambda}$ attached to the Kronecker coefficients are quasipolynomials \citep{GCT6}. 
This means that for any fixed $\lambda$, $\mu$, $\nu$ there exist an integer $k$ and polynomials $p_1$, $p_2$, \ldots, $p_k$ such that for any $N \geq 1$, $\widetilde{g}_{\mu,\nu}^{\lambda}(N)=p_i(N)$ when $N \equiv i \mod k$. 
Mulmuley's SH conjecture stated that for any such description, $g_{\mu,\nu}^{\lambda}=0 \Leftrightarrow F_1=0$. 
The rightmost system of conditions in \eqref{eq:conditions} above provides a family of counterexamples to this conjecture \citep{Briand:Orellana:Rosas:SH}. The discovery of these counterexamples led \citet{Mulmuley:erratum} to propose a weaker form of the conjecture SH, still strong enough for the aims of Geometric Complexity Theory.


\bibliographystyle{plainnat}
\def\cprime{$'$}

\end{document}